\documentclass[12pt]{article}

\usepackage[T1]{fontenc}
\usepackage[utf8]{inputenc}
\usepackage{amsmath,amssymb,amsthm,mathtools}
\usepackage[hidelinks]{hyperref}

\newtheorem{theorem}{Theorem}

\DeclareMathOperator{\diag}{diag}

\newcommand{\R}{\mathbb{R}}
\newcommand{\C}{\mathbb{C}}

\newcommand{\beq}{\begin{equation}}
\newcommand{\eeq}{\end{equation}}

\newcommand{\entrygeq}{\geq}

\title{On the positive semidefinteness of a class of Hermitian
Cauchy-like matrices}
\author{Augusto Ferrante\thanks{Dipartimento di Ingegneria dell'Informazione, Universit\`a di Padova, Via Gradenigo 6/B, 35131 Padova. email: {\texttt{augusto.ferrante@unipd.it}}}}
\date{}

\begin{document}

\maketitle

\begin{abstract}
We establish the positive semidefiniteness of a class of Hermitian Cauchy-like matrices associated with real Hurwitz polynomials having distinct zeros. Writing the zeros as $-\lambda_1,\ldots,-\lambda_n$, the matrix entries are defined through ratios of elementary symmetric polynomials in the variables $\lambda_i$. The result  extends an earlier positivity theorem obtained under the assumption that all zeros are real. We first show that the coefficients appearing in the denominators are nonzero, so that the matrices are well defined. We then construct an explicit congruence between each matrix and the solution of a Lyapunov equation whose state matrix is in companion form. Positive semidefiniteness follows from a recent result on such equations with entrywise nonnegative right-hand-side data. This approach establishes the extension to complex zeros through the connection between structured matrices and Lyapunov equations.
\end{abstract}

\section{Introduction}

Let $\sigma(s)$ 
be a  Hurwitz monic polynomial of degree $n$ with real coefficients and assume that
$\sigma(s)$ has $n$ distinct zeros $-\lambda_1,-\lambda_2, \dots,-\lambda_n\in\C$.
Since $\sigma(s)$ is Hurwitz, each $\lambda_i$ has positive real part.
For $m=0,1,\dots, n$, denote by $\sigma_{m}$ the \emph{elementary symmetric polynomial of degree $m$}
in the variables $\{\lambda_1, \dots, \lambda_n\}$, i.e. 
$\sigma_0:=1$ and $\sigma_m
:=
\sum_{1\le i_1<i_2<\cdots<i_m\le n}
\lambda_{i_1}\lambda_{i_2}\cdots \lambda_{i_m}$, for $m=1,\dots, n$.
Then we have
\beq\label{def-sigma}
\sigma(s)=s^n+\sigma_1s^{n-1}+\dots+ \sigma_{n-1}s + \sigma_n.
\eeq
For $m=0,1,\dots, n-1$, denote $\sigma_{m}^{(j)}$ the  elementary symmetric polynomial of degree $m$
in the variables $\{\lambda_1, \dots, \lambda_n\} \setminus \{\lambda_j\}$.
For $l,k\in\{1,2,\dots,n\}$, define the  complex numbers
\beq\label{rikl}
r_i(k,l):=\frac{\sigma_{n-k}^{(i)}}{\sigma_{n-l}^{(i)}}
\eeq
and consider the class of Hermitian Cauchy-like matrices  $C(k,l)\in\C^{n\times n}$
whose elements are 
\beq\label{defcij}
C_{ij}(k,l)=\frac{r_i(k,l)+\overline{r_j(k,l)}}{\lambda_i+\overline{\lambda_j}},\quad i,j=1,\dots,n.
\eeq

Our main result is the following theorem:

\begin{theorem}\label{main-res}
Let $\sigma(s)$ given by (\ref{def-sigma})
be a real Hurwitz  polynomial of degree $n$ having $n$ distinct zeros $-\lambda_1,-\lambda_2, \dots,-\lambda_n\in\C$.
Then the matrix $C(k,l)\in\C^{n\times n}$ whose elements are defined by (\ref{defcij}) is positive semidefinite for all $l,k\in\{1,2,\dots,n\}$.
\end{theorem}

This theorem generalizes the main result of a recent paper \cite{ferrante2026positivityclasscauchylikematrices} in which the same result has been proven under the further assumption that all the $\lambda_i$ are real.
The proof in that case was based on a direct and detailed analysis of the matrix.
It appears, however, that this analysis cannot be generalized to the complex case.
On the other hand, a very recent result, \cite{ferrante2026dualityreformulationcompanionmatrixlyapunov}, allows to prove Theorem \ref{main-res}
via an indirect route through the connection with a property of Lyapunov equations
with state matrices in companion form.
Matrices of form $C(k,l)$
whose elements are defined in 
\eqref{defcij} can be considered to be a variant of Cauchy
matrices \cite[Chap. 3]{Pan2001}.
More specifically our matrices are related to a very well studied class of matrices:
the so-called Kwong (or anti-Loewner) matrices.
The latter are matrices $W$ whose elements have the form
$$
W_{ij}=\frac{f(\lambda_i)+f(\lambda_j)}{\lambda_i+\lambda_j},
$$
where $\lambda_i$ are distinct real numbers and $f$ is a real-valued function.
The matrices studied in this paper address the more general case of complex
Hermitian matrices (instead of real symmetric matrices) and can be viewed as dual to Kwong matrices.
In fact, setting $\Lambda^{(i)}:= \{\lambda_j: \ j=1,\dots,n, \ j\neq i\}$,
we can write 
$$C_{ij}=\frac{f(\Lambda^{(i)})+\overline{f(\Lambda^{(j)})}}{\lambda_i+\overline{\lambda_j}}$$
where $f$ is a complex-valued function. 

As a variation on the classical 
results of Loewner
\cite{Loewner1934}, a quite extensive literature has analyzed Kwong matrices and the properties of $f$ guaranteeing their positivity \cite{Kwong1989,BhatiaSano2009,Audenaert2011,Morishita2014,BhatiaJain2023}.
On the other hand, the matrices of the form discussed here have been addressed only in
\cite{ferrante2026positivityclasscauchylikematrices}. This paper may be viewed as a second step toward filling this gap.

{\bf Notation.}
For $X\in\R^{n\times m}$, the notation \(X\entrygeq0\) (\(X> 0\)) means that \(X\) is entrywise nonnegative (entrywise positive).
For a complex matrix $P\in\C^{n\times m}$, $\overline{P}$, $P^\top$ and $P^*$
denote, respectively, the complex conjugate, the transpose and  the conjugate transpose of $P$.
If $X=X^*\in\C^{n\times n}$, \(X\succeq 0\) (\(X\succ 0\)) denotes positive semidefiniteness (positive definiteness).
We denote by $\mathrm{i}$ the imaginary unit.

\section{Proof of Theorem \ref{main-res}}
First of all, we need to prove that $C(k,l)$ is well defined
or, equivalently,  that for all $l=1,\dots,n$,
\begin{equation}\label{eq:denominators-nonzero}
\sigma_{n-l}^{(i)}\ne0.
\end{equation}
To this aim, define the monic Hurwitz polynomials
\begin{equation}\label{eq:qi}
q_i(s):=\prod_{\small{
\begin{array}{c}j=1\\[-1mm]
j\neq i\end{array}}}^{n}(s+\lambda_j)=\frac{\sigma(s)}{s+\lambda_i},\qquad i=1,\dots,n.
\end{equation}
By definition, the coefficients of $q_i$ are precisely the complex numbers $\sigma_m^{(i)}$.

If $\lambda_i$ is real, then $q_i(s)$ is a real Hurwitz polynomial, so all its
coefficients are strictly positive.

If $\lambda_i=\alpha+\mathrm{i}\beta$, with $\beta\ne0$, write
\begin{equation}\label{eq:complex-factorization}
\sigma(s)=(s+\lambda_i)(s+\overline{\lambda_i})b_i(s),
\end{equation}
where $b_i(s)$ is a real polynomial of degree equal to $n-2$; moreover,
$b_i(s)$ is monic and Hurwitz  and therefore all its coefficients are strictly positive. Then
\begin{equation}\label{eq:qi-complex}
q_i(s)=(s+\overline{\lambda_i})b_i(s)=sb_i(s)+\alpha b_i(s)-\mathrm{i}\beta b_i(s).
\end{equation}
Thus, the leading coefficient of $q_i(s)$ is $1\neq 0$. The other coefficients have imaginary part equal to $-\beta$ multiplied by a strictly positive
coefficient of $b_i(s)$ and hence are also nonzero.
Hence no coefficient of $q_i$ vanishes. Therefore $\sigma_m^{(i)}\neq 0$
for all $m=0,\dots, n-1$ and $i=1,\dots,n$.
In conclusion, $C(k,l)$ is well defined.
To show that it is positive semidefinite,
define
\begin{equation}\label{comp-matrix}
A:=
\begin{bmatrix}
0&1&0&\cdots&0\\
0&0&1&\ddots&\vdots\\
\vdots&&\ddots&\ddots&0\\
0&\cdots&0&0&1\\
-\sigma_n&-\sigma_{n-1}&\cdots&-\sigma_2&-\sigma_1
\end{bmatrix}\in \R^{n\times n}
\end{equation}
to be the companion matrix associated with the real  Hurwitz polynomial
$\sigma(s)$.

Consider the Lyapunov equation
\begin{equation}\label{eq:lyap-primal}
  AP+PA^{\top}=-Q
\end{equation}
where $Q=Q^\top\in\R^{n\times n}$.
As shown in \cite{ferrante2026dualityreformulationcompanionmatrixlyapunov}, if $Q$ is entrywise nonnegative then the solution  $P$ 
of \eqref{eq:lyap-primal} is positive semidefinite.
In particular, select
\beq\label{elementary-Q}
Q=
Q(k,l):=\varepsilon_k \varepsilon_l^\top+\varepsilon_l \varepsilon_k^\top,
\eeq
 where $\varepsilon_k$ denotes $k$-th canonical vector and $k,l=1,\dots,n$.
 Then $Q(k,l)=Q(k,l)^\top$ is entrywise nonnegative.
 Hence, the associated solution $P(k,l)$ of \eqref{eq:lyap-primal}, which is given by, 
 \begin{equation}\label{eq:Pkl}
P(k,l)=\int_0^\infty e^{At}Q_{kl}e^{A^\top t}\,dt
\end{equation}
is positive semidefinite: 
\beq\label{P-psd}
P(k,l)\succeq 0,\qquad \forall k,l=1,\dots n.
\eeq
It is well-known and easy to check that, for $i=1,\dots,n$, the 
 left eigenvector of $A$ associated with the eigenvalue $-\lambda_i$ is
\begin{equation}\label{eq:vi}
v_i^\top:=
\begin{bmatrix}
\sigma_{n-1}^{(i)} & \sigma_{n-2}^{(i)} & \cdots & \sigma_0^{(i)}
\end{bmatrix}^{\!\top}.
\end{equation}
By taking \eqref{elementary-Q} into account, we see that
\begin{equation}\label{eq:eigenvector-Qkl}
v_i^\top Q(k,l) \overline{v_j}
=
\sigma_{n-k}^{(i)}\overline{\sigma_{n-l}^{(j)}}
+\sigma_{n-l}^{(i)}\overline{\sigma_{n-k}^{(j)}}.
\end{equation}
As a consequence, 
\begin{equation}\label{eq:eigenvector-Pkl}
\begin{aligned}
v_i^\top P(k,l)\overline{v_j}
&=
\int_0^\infty v_i^\top e^{At}Q(k,l)e^{A^\top t}\overline{v_j}\,dt\\
&=\int_0^\infty e^{-\lambda_it}v_i^\top Q(k,l)\overline{v_j}e^{-\overline{\lambda_j} t}\,dt\\
&=\int_0^\infty e^{-\lambda_it}e^{-\overline{\lambda_j} t}\,dt
\left[\sigma_{n-k}^{(i)}\overline{\sigma_{n-l}^{(j)}}
+\sigma_{n-l}^{(i)}\overline{\sigma_{n-k}^{(j)}}\right]\\
&=\frac{
\sigma_{n-k}^{(i)}\overline{\sigma_{n-l}^{(j)}}
+\sigma_{n-l}^{(i)}\overline{\sigma_{n-k}^{(j)}}
}{\lambda_i+\overline{\lambda_j}}.
\end{aligned}
\end{equation}
Set
\begin{equation}\label{eq:V}
V:=
\begin{bmatrix}
\overline{v_1}&\cdots&\overline{v_n}
\end{bmatrix}
\end{equation}
and observe that the last member of (\ref{eq:eigenvector-Pkl}) is equal to element in row $i$ and column $j$ of $V^*P(k,l) V$.
Consequently, we have:
$$
C(k,l)=D^{-1}V^*P(k,l) V D^{-*}$$
where
\begin{equation}\label{eq:D}
D:=\diag\!\left(
\sigma_{n-l}^{(1)},\dots,\sigma_{n-l}^{(n)}
\right).
\end{equation}
Hence,
\begin{equation}\label{eq:C-P-equivalence}
C(k,l)\succeq 0
\quad\Longleftrightarrow\quad
P(k,l)\succeq 0,
\end{equation}
because $D$ and $V$ are nonsingular. This concludes the proof in view of \eqref{P-psd}.
\qed

\section{Conclusions}
We have proved that the Hermitian Cauchy-like matrices $C(k,l)$ associated with a real Hurwitz polynomial with distinct zeros are positive semidefinite for every $k,l\in{1,\ldots,n}$. This extends the previously established result for real zeros to the complex setting permitted by real polynomial coefficients.

The proof rests on two observations. First, the coefficients of the polynomials obtained by removing one linear factor are nonzero, ensuring that all defining ratios exist. Second, each $C(k,l)$ is congruent to the solution of a companion-matrix Lyapunov equation with a symmetric, entrywise nonnegative forcing matrix. The congruence therefore transfers the known positivity property of these Lyapunov solutions to the matrices considered here. Beyond proving the extension, this representation makes explicit the connection between the elementary symmetric polynomial construction and companion-matrix dynamics, and broadens the class of positive semidefinite structured matrices related to Kwong matrices.

\end{document}